\documentclass[12pt]{amsart}
\usepackage{amssymb, amstext, amscd, amsmath}
\makeatletter
\def\@cite#1#2{{\m@th\upshape\bfseries%
[{#1\if@tempswa{\m@th\upshape\mdseries, #2}\fi}]}} \makeatother
\theoremstyle{plain}
\newtheorem{thm}{Theorem}[section]
\newtheorem{cor}[thm]{Corollary}
\newtheorem{prop}[thm]{Proposition}

\theoremstyle{definition}
\newtheorem{rem}[thm]{Remark}
\newtheorem{defn}[thm]{Definition}

\newtheorem{eg}[thm]{Example}
\newcommand{\Prf}{\noindent\textbf{Proof.\ }}
\newcommand{\bx}{\hfill$\blacksquare$\medbreak}

\newcommand{\ca}{\mathrm{C}^*}

\newcommand{\ol}{\overline}
\newcommand{\td}{\widetilde}

\newcommand{\wot}{\textsc{wot}}

\newcommand{\bbC}{{\mathbb{C}}}
\newcommand{\bbF}{{\mathbb{F}}}

\newcommand{\bbN}{{\mathbb{N}}}

\newcommand{\bbT}{{\mathbb{T}}}
\newcommand{\bbZ}{{\mathbb{Z}}}

  \newcommand{\A}{{\mathcal{A}}}
  \newcommand{\B}{{\mathcal{B}}}

  \newcommand{\E}{{\mathcal{E}}}

\renewcommand{\H}{{\mathcal{H}}}
  
  \newcommand{\J}{{\mathcal{J}}}
  \newcommand{\K}{{\mathcal{K}}}

  \newcommand{\M}{{\mathcal{M}}}

\renewcommand{\P}{{\mathcal{P}}}

  \newcommand{\V}{{\mathcal{V}}}
  \newcommand{\W}{{\mathcal{W}}}

\newcommand{\fL}{{\mathfrak{L}}}

\newcommand{\fR}{{\mathfrak{R}}}

\newcommand{\qand}{\quad\text{and}\quad}

\newcommand{\qfor}{\quad\text{for}\quad}
\newcommand{\qforal}{\quad\text{for all}\quad}

\newcommand{\ran}{\operatorname{Ran}}
\newcommand{\rank}{\operatorname{rank}}
\newcommand{\spn}{\operatorname{span}}

\newcommand{\fgeeplus}{\bbF^+\!(G)}

\newcommand{\Lg}{{\mathfrak{L}}_G}
\newcommand{\Rg}{{\mathfrak{R}}_G}

\newcommand{\bofh}{\B(\H)}
\newcommand{\rowt}{(T_e)_{e\in E(G)}}
\newcommand{\rows}{(S_e)_{e\in E(G)}}
\newcommand{\rowl}{(L_e)_{e\in E(G)}}

\newcommand{\flgee}{\fL_G}
\newcommand{\frgee}{\fR_G}

\begin{document}

\title[Wold Decomposition]{Applications of the Wold decomposition to the study of
row contractions associated with directed graphs}
%
%
\author[E. Katsoulis]{Elias~Katsoulis}
\address{Department of Mathematics\\East Carolina University\\
Greenville, NC 27858\\USA} \email{KatsoulisE@mail.ecu.edu}

\author[D.W. Kribs]{David~W.~Kribs}
\address{Department of Mathematics and Statistics\\University of Guelph\\
Guelph, Ontario\\CANADA N1G 2W1} \email{dkribs@uoguelph.ca}
\begin{abstract}
Based on a Wold decomposition for  families of partial isometries
and projections of Cuntz-Krieger-Toeplitz-type, we extend several
fundamental theorems from the case of single vertex graphs to the
general case of countable directed graphs with no sinks. We prove
a Szego-type factorization theorem for CKT families, which leads
to information on the structure of the unit ball in free
semigroupoid algebras, and show that joint similarity implies
joint unitary equivalence for such families. For each graph we
prove a generalization of von Neumann's inequality which applies
to  row contractions of operators on Hilbert space which are
related to the graph in a natural way. This yields a functional
calculus determined by quiver algebras and free semigroupoid
algebras. We establish a generalization of Coburn's theorem for
the $\ca$-algebra of a CKT family, and prove a universality
theorem for $\ca$-algebras generated by these families. In both
cases, the $\ca$-algebras generated by quiver algebras play the
universal role.
\end{abstract}

\date{}
\maketitle

\section{Introduction}\label{S:intro}

In \cite{JK1}, the second author and Jury derived a Wold
decomposition for families of partial isometries which satisfy the
`Cuntz-Krieger-Toeplitz' directed graph relations, what we call
the  $(\dagger)$ relations. This theorem may also be obtained as a
special case of the Wold decomposition recently established by
Muhly and Solel \cite{MS1} for the more general setting of induced
representations of tensor algebras. In \cite{JK1} it was used to
investigate the internal structure of free semigroupoid algebras
\cite{JP,JK1,KK,KP1,KP2} and a dilation theory was built around it
in \cite{JK2,MS2}. In this paper, we apply this Wold decomposition
to establish extensions of a number of fundamental theorems to the
case of general countable directed graphs with no sinks, whereby
the extended results can be regarded as occurring for single
vertex graphs.

The next section ($\S 2$) contains background information.  In $\S
3$ we include an expanded discussion of the Wold decomposition and
extend a folklore result for isometries on Hilbert space: We prove
that CKT families $\{ S_e, P_x\}$ of partial isometries and
projections are jointly similar precisely when they are jointly
unitarily equivalent. We also establish a Szego-type factorization
theorem \cite{Rov,Pop4} for CKT families. In $\S 4$ we use the
Szego theorem to obtain detailed information \cite{DKP} on the
unit ball of semigroupoid algebras which are partly free
\cite{KP1,KP2}.

Given a countable directed graph $G$ with no sinks, in $\S 5$ we
show there is a von Neumann inequality for the set of  row
contractions $T = (T_1,\ldots,T_n)$ of operators on Hilbert space
such that the joint actions of the $T_i$ are related to $G$ in a
natural way. In the case of single vertex graphs, we recover the
original form of the von Neumann inequality \cite{vN,SzF} (when
$n=1$) and Popescu's version \cite{Popvn} (when $2\leq n \leq
\infty$). If $n\geq 2$,  the estimates obtained  are sharper than
the estimates of \cite{Popvn} in general. In fact, we show these
estimates respect the partial ordering on directed graphs induced
by graph deformations. An immediate consequence of the von Neumann
inequality  is the existence of functional calculus ($\S 6$) for
row contractions based on quiver algebras
\cite{KK,KP1,KP2,Mu,MS1,MS2} and free semigroupoid algebras. There
is a long history of both positive and negative results associated
with von Neumann inequalities for the multivariable setting, and
this work fits in with such efforts. See
\cite{Ando,Arv,BC,Boz,CD,Drury,Hol,MS2,Popvn} for some of the
different perspectives. In connection with the current work, see
\cite{MS2} where a more abstract form of the von Neumann
inequality is obtained.

In the final section, we apply the Wold decomposition to obtain a
pair of universality theorems for directed graph $\ca$-algebras.
In particular, a generalization of Coburn's theorem
\cite{Cob,Pop2} for isometries on Hilbert space is established:
The $\ca$-algebra generated by any CKT family which satisfies a
natural non-unitary condition is isomorphic to the $\ca$-algebra
generated by the corresponding quiver algebra.  Further, we prove
that the $\ca$-algebras generated by quiver algebras are the
universal $\ca$-algebras for the CKT relations.

\section{Preliminaries}\label{S:prelim}

Let $G$ be a countable directed graph with directed edges $\E(G)$
and vertices $\V(G)$. Given a finite path $w$ in $G$, we write
$w=ywx$ when the source vertex of $w$ is $s(w) =x$ and the range
vertex of $w$ is $r(w)=y$. (For vertices $x\in \V(G)$ put
$r(x)=x=s(x)$.) Say that $G$ has a {\it source} at $x\in \V(G)$ if
there are no edges which finish at $x$; that is, $r(e)\neq x$ for
$e\in\E(G)$. Similarly, $G$ has a {\it sink} at $x\in \V(G)$ if
$s(e)\neq x$ for $e\in \E(G)$. Define the {\it free semigroupoid}
$\fgeeplus$ to be the set of all vertices $\V(G)$ and all finite
paths in $G$, written in reduced form, with the natural operations
of concatenation of paths allowed by the graph structure. A path
$w = e_k \cdots e_1 \in \fgeeplus$ is said to be a \textit{loop}
if $s(w) = r(w)$. If in addition, $r(e_i) \neq r(e_j)$, for $i
\neq j$, then $w$ is said to be \textit{vertex-simple}. A
vertex-simple loop $w = e_k  \cdots e_1$, $e_i \in \E(G)$, has an
{\it entrance} if there exists an edge $f \in \E(G)$ such that
$r(f) = r(e_i)$, for some $i$, but $f \neq e_i$. Finally, let
$\P_G^+$ be the {\it path algebra} generated by $G$; this is the
set of polynomials with complex coefficients and monomials
belonging to $\fgeeplus$.

Given a family of non-zero partial isometries $\{S_e :
e\in\E(G)\}$ and projections $\{P_x : x\in\V(G)\}$ which act on a
common Hilbert space, consider the following relations:
\[
(\dagger)  \left\{
\begin{array}{lll}
(i)  & P_x P_y = 0 & \mbox{for all $x,y \in \V(G)$, $ x \neq y$}  \\
(ii) & S_{e}^{*}S_f = 0 & \mbox{for all $ e, f \in \E(G)$, $e \neq f $}  \\
(iii) & S_{e}^{*}S_e = P_{s(e)} & \mbox{for all $e \in \E(G)$}      \\
(iv)  & \sum_{r(e)=x}\, S_e S_{e}^{*} \leq P_{x} & \mbox{for all
$x \in \V(G)$}
\end{array}
\right.
\]
If equality is achieved in $(iv)$ for all $x\in\V(G)$, then we
refer to the relations as $(\ddagger)$. Observe that $(\ddagger)$
can only occur for $G$ with no sources. We mention that the
operator algebras generated by families $\{S_e, P_x\}$ which
satisfy $(\dagger)$ have been the focus of intense recent
interest. They include so-called Cuntz-Krieger directed graph
$\ca$-algebras, free semigroupoid algebras, quiver algebras, etc.
(The references
\cite{BHRS,CrKum,JP,JK1,KK,KumP,KPR,KumP2,KP1,KP2,Mu,MS1,MS2,Sol,Szy}
give a starting point for the interested reader.) Let us mention
at this point the existence of the universal $\ca$-algebra
$\ca(G)$ \cite{KumP,KPR,KumP2} for the relations $(\ddagger)$ for
a given directed graph $G$. From these operator algebra
motivations, we refer to $(\ddagger)$ as the Cuntz-Krieger (CK)
relations, and $(\dagger)$ as the Cuntz-Krieger-Toeplitz relations
(CKT).

We next describe the prototypical `pure model' \cite{JK2,JK1} in
the Wold decomposition used here. (More details are given in the
next section.) Let $G$ be  a countable directed graph (with no
sinks) and define $\H_G$ to be the Hilbert space with orthonormal
basis given by $\{\xi_w : w\in\fgeeplus\}$. This generalized `Fock
space' associated with a directed graph was introduced by Muhly
\cite{Mu}. For $e\in \E(G)$ define partial isometries $L_e$ on
$\H_G$ by
\begin{eqnarray}\label{ledefn}
L_e \xi_w = \left\{ \begin{array}{cl} \xi_{ew} & \mbox{if $r(w) = s(e)$} \\
0 & \mbox{otherwise}
\end{array}\right.
\end{eqnarray}
Observe that $L_e^*L_e = P_{s(e)}$ is the projection onto
$\spn\{\xi_w : r(w) = s(e)\}$, and that the family $\{L_e, P_x\}$
satisfies $(\dagger)$. For $x\in\V(G)$, we define the `tree
component' subspace $\H_x = \spn\{ \xi_w : w=wx\}$. Notice that
these subspaces are reducing for the $L_e$. In fact, the subspaces
$\{ \H_x : x\in\V(G)\}$ form the {\it unique} family of minimal
non-zero joint reducing subspaces for the $L_e$ \cite{KP1}.

Given a directed graph $G$, the associated {\it quiver algebra}
$\A_G$ is the norm-closed algebra generated by $\{ L_e, P_x :
e\in\E(G), x\in\V(G)\}$, where the $P_x$ here are the vertex
projections associated with the pure model (\ref{ledefn}). The
associated {\it free semigroupoid algebra} $\flgee$ is the weak
operator topology closure of $\A_G$. The commutant of $\flgee$
coincides with $\flgee^\prime = \frgee$ \cite{KP1}, the
$\wot$-closed algebra generated by partial isometries $R_e$,
defined on $\H_G$ (as in (\ref{ledefn})) by $R_e \xi_w =\xi_{we}$,
together with their initial projections which we denote by $Q_x$,
$x\in\V(G)$. The projections $Q_x$ are precisely the projections
onto the tree component subspaces $\H_x$.

We mention that the graph $G$ was proved to be a complete
invariant for $\A_G$ (and $\flgee$) up to unitary equivalence by
Kribs and Power \cite{KP1}. More recently, Solel \cite{Sol} proved
$G$ is an invariant for isometric isomorphism and, independently,
the authors \cite{KK} proved that $G$ is a banach algebra
isomorphism invariant for $\A_G$. (In fact, when $G$ has no
sources or no sinks it was proved in \cite{KK} that $G$ is an
algebraic isomorphism invariant for $\A_G$.) This is a true
departure from the $\ca$-algebra case, where, for instance, it is
not hard to see that $G$ is not a $\ast$-isomorphism invariant for
$\ca(G)$. (e.g. For many non-isomorphic graphs, $\ca(G)$ is
$\ast$-isomorphic to the set of compact operators.)

\section{Wold Decomposition}\label{S:wold}

The Wold decomposition from \cite{JK1} (which assumed $G$ had no
sinks) asserts that every set of partial isometries $\{ S_e :
e\in\E(G)\}$ which satisfies $(\dagger)$ on a Hilbert space $\H$
is jointly unitarily equivalent to a direct sum
\begin{eqnarray}\label{wold}
\,\,\,\,\,\,\,\,\, S_e \simeq
 V_e \oplus \Big( \sum_{x\in V(G)}\oplus L_e^{(\alpha_x)}\Big|_{\H_x^{(\alpha_x)}}
\Big) \qfor e\in E(G),
\end{eqnarray}
where  $V = (V_e)_{e\in\E(G)}$ satisfies $(\ddagger)$, and hence
determines a representation of a Cuntz-Krieger directed graph
$\ca$-algebra. As indicated above, this result can be derived from
a special case of the Wold decomposition from \cite{MS1}. The
restricted ampliations of the $L_e$ in this decomposition are said
to form the {\it pure part} of the dilation. (Notice that the pure
model (\ref{ledefn}) is captured with $V_e \equiv 0$ and $\alpha_x
\equiv 1$.) The $\alpha_x$ are called the {\it vertex
multiplicities} in the dilation and are computed from
$\{S_e,P_x\}$ via the equations
\begin{eqnarray}\label{vertexmult}
\alpha_x = \rank \Big( P_x \big( I - \sum_e S_e S_e^* \big) \Big)
\qfor x\in\V(G).
\end{eqnarray}

More specifically, the joint unitary equivalence (\ref{wold})
arises from the spatial decomposition $\H = \H_c \oplus \H_p$
where $\H_c = \H_p^\perp$  and $\H_p =\sum_{w\in\fgeeplus}\oplus
w(S) \W$ are joint reducing subspaces for the $S_e$, with $\W$ the
wandering subspace for $S$ given by
\[
\W = \ran \big( I - \sum_e S_e S_e^*\big) = \bigcap_e \ker S_e^*.
\]
Observe for $x\in\V(G)$ that $P_x\W = \ran \big( P_x -
\sum_{r(e)=x} S_eS_e^*\big)$ is a subspace of $\W$. The subspaces
$\H_c$ and $\H_p$ have the alternate descriptions
\begin{eqnarray*}
\H_p = \Big\{ \xi\in\H\,\,\,:\,\,\, \lim_{d\rightarrow \infty}
\sum_{w\in\fgeeplus_d} ||w(S)^* \xi||^2 =0 \Big\}
\end{eqnarray*}
\begin{eqnarray}\label{ckpart}
\H_c = \bigcap_{d\geq 1} \Big\{ \xi\in\H: \sum_{w\in\fgeeplus_d}
w(S)w(S)^* \xi=\xi \Big\},
\end{eqnarray}
where $\bbF^+(G)_d$ is the set of paths inside $\fgeeplus$ of
length $d$.

As a direct consequence of the Wold decomposition, we prove a
generalization of the well-known fact that similar isometries are
unitarily equivalent.  This result also generalizes Popescu's
theorem \cite{Pop2} for the case of $n$-tuples of isometries with
orthogonal ranges, which occurs below when $G$ is a single vertex
graph with $n\geq 2$ edges.

\begin{thm}\label{simunit}
Let $G$ be a countable directed graph with no sinks. If
$\{S_e,P_x\}$ and $\{S^\prime_e, P_x^\prime\}$ both satisfy
$(\dagger)$ for $G$ and are jointly similar, then they are jointly
unitarily equivalent. That is, if there is an invertible operator
$A$ with
\begin{eqnarray*}
A S_e A^{-1} = S_e^\prime & \qfor e\in\E(G) \\
A P_x A^{-1} = P_x^\prime & \qfor x\in\V(G),
\end{eqnarray*}
then there is a unitary operator $U$ such that
\begin{eqnarray*}
U S_e U^* = S_e^\prime & \qfor e\in\E(G) \\
U P_x U^* = P_x^\prime & \qfor x\in\V(G).
\end{eqnarray*}
\end{thm}

\Prf Without loss of generality assume both families of operators
act on $\H$. Let $\H = \H_c \oplus \H_p = \H_c^\prime \oplus
\H_p^\prime$ be the spatial decompositions associated with the
Wold decompositions (\ref{wold}) for the $S_e$ and $S_e^\prime$
respectively.

For $e\in\E(G)$ let
\[
S_e|_{\H_c} \equiv V_e, \quad S_e|_{\H_p} \equiv W_e, \quad
S_e^\prime|_{\H_c^\prime} \equiv V_e^\prime, \quad
S_e^\prime|_{\H_p^\prime} \equiv W_e^\prime.
\]
By hypothesis, for $e\in\E(G)$ we have
\[
\ker S_e^* = \ker(A^* S_e^{\prime *} (A^*)^{-1}) = A^* \big(\ker
(S_e^{\prime *})\big).
\]
Let $\W$ and $\W^\prime$ be the wandering subspaces for $S$ and
$S^\prime$ respectively. Then
\[
A^* (\W^\prime) =A^* \Big( \bigcap_e \ker S_e^{\prime *} \Big) =
\bigcap_e \ker S_e^* = \W.
\]
In particular, $\dim\W = \dim \W^\prime$. But more is true. Recall
that the vertex multiplicities in the pure parts of the Wold
decompositions for $S$ and $S^\prime$ are given by $\alpha_x=\dim
P_x \W$ and $\alpha_x^\prime = \dim P_x \W^\prime$. By hypothesis
and the previous identity we have
\[
P_x \W = P_x A^* \W^\prime = A^* P_x^\prime \W^\prime \qfor
x\in\V(G),
\]
and hence the invertibility of $A$ yields $\alpha_x =
\alpha_x^\prime$ for $x\in\V(G)$.

Thus, it is easily seen from the Wold decomposition that there is
a unitary $U_1:\H_p \rightarrow \H_p^\prime$ which intertwines the
pure parts of $S$ and $S^\prime$;
\begin{eqnarray}\label{commute0}
U_1 W_e U_1^* = W_e^\prime  \qfor e\in\E(G).
\end{eqnarray}
Since $G$ has no sinks, every vertex projection $P_x$ is the
initial projection for some $S_e$. Further, recall that $\H_p$
(respectively $\H_p^\prime$) is a reducing subspace for the $S_e$
(respectively $S_e^\prime$). It follows that the restricted vertex
projections $P_x|_{\H_p}$ and $P_x^\prime|_{\H_p^\prime}$ are also
intertwined by $U_1$.

We now show that the restrictions to $\H_c, \H_c^\prime$ may be
spatially intertwined as well. As $AS_e = S_e^\prime A$ and $S_e
A^{-1} = A^{-1} S_e^\prime$ for $e\in \E(G)$, equation
(\ref{ckpart}) applied for the $S_e$ and the $S_e^\prime$ implies
that $A \H_c = \H_c^\prime$. Let $A_0 = A|_{\H_c}:
\H_c\rightarrow\H_c^\prime$ and observe that $A_0$ is invertible.
Then we have
\[
A_0 V_e = V_e^\prime A_0 \qfor e\in\E(G).
\]
In particular, equation (\ref{ckpart}) gives
\[
A_0 = A_0 (I_{\H_c}) = A_0 \Big( \sum_e V_e V_e^* \Big) = \sum_e
V_e^{\prime} A_0 V_e^*.
\]
Thus the $(\dagger)$ relations and the hypothesis implies that for
$e\in \E(G)$,
\[
V_e^{\prime *} A_0 = V_e^{\prime *} V_e^\prime A_0 V_e^* =
P_{s(e)}^\prime A_0 V_e^* = A_0 P_{s(e)} V_e^* = A_0 V_e^*.
\]
It follows from the adjoint of this identity that
\begin{eqnarray}\label{commute}
V_e A_0^* A_0 = A_0^* V_e^\prime A_0 = A_0^* A_0 V_e \qfor
e\in\E(G).
\end{eqnarray}

Consider the polar decomposition of $A_0= U_2 B$, where $U_2: \H_c
\rightarrow \H_c^\prime$ is unitary and $B = (A_0^* A_0)^{1 / 2}$
is an invertible operator on $\H_c$. By (\ref{commute}) we have
\[
V_e B = B V_e \qfor e\in\E(G),
\]
and it follows that
\begin{eqnarray*}
V_e^\prime U_2 = V_e^\prime U_2 B B^{-1} &=& V_e^\prime A_0
B^{-1}\\ &=& A_0 V_e B^{-1} = U_2 B V_e B^{-1} = U_2 V_e.
\end{eqnarray*}
Thus
\begin{eqnarray}\label{commute1}
U_2 V_e = V_e^\prime U_2 \qfor e\in\E(G).
\end{eqnarray}

Therefore, we may define a unitary $U = U_1 \oplus U_2$ on $\H$ so
that (\ref{commute0}) and (\ref{commute1}) imply
\[
U S_e = S_e^\prime U^* \qfor e\in\E(G).
\]
As $G$ has no sinks, it follows that the projections $P_x$ and
$P_x^\prime$ are also intertwined by the unitary $U$. \bx


We now present another application of the Wold decomposition that
allows us to prove a Szego-type factorization theorem. This
theorem generalizes earlier results of Rosenblum and Rovnyak
\cite[Theorem 3.4]{Rov} and Popescu \cite{Pop4}.

\begin{thm}
Let $G$ be a countable directed graph with no sinks and let $\{
P_x\}_{x \in \V(G)}$ and $\{ S_e \}_{e \in \E(G)}$ be sets of
projections and partial isometries which satisfy $(\dagger)$.
Assume that $\{ S_e \}_{e \in \E(G)}$ is pure. If $ Y \in B(\H )$
is a positive invertible operator then the following assertions
are equivalent:
\begin{itemize}
\item[(i)] $Y = A^{*}A$ for some $A \in B(\H)$ which commutes with
$\{ P_x\}_{x \in \V(G)}$ and $\{ S_e \}_{e \in \E(G)}$.
\item[(ii)] $S_{e}^{*}YS_f = YS_{e}^{*} S_f = S_{e}^{*} S_f Y$ for
all $e,f \in \E(G)$.
\end{itemize}
\end{thm}

\Prf We  need to prove that (ii) implies (i). Let $T_u \equiv
Y^{1/2}S_u Y^{-1/2}$ for $u \in \fgeeplus$. We will show that the
family $\{ T_e \}_{e \in \E(G)}$ is jointly unitarily equivalent
to $\{ S_e \}_{e \in \E(G)}$. Towards this end notice that
$T_{e}^{*}T_e = P_{s(e)}$ and $T_e T_{e}^{*} \leq P_{r(e)}$.
Indeed,
\begin{eqnarray*}
T_{e}^{*}T_e  &=&  Y^{-1/2}S_{e}^{*}YS_eY^{-1/2} \\
              &=&    Y^{-1/2}Y Y^{-1/2} P_{s(e)}
              = P_{s(e)} = S_e^* S_e,
\end{eqnarray*}
and a similar calculation shows the claimed inequality. Therefore,
the collections $\{ P_x\}_{x \in \V(G)}$ and $\{ T_e \}_{e \in
\E(G)}$ satisfy the requirements of the Wold Decomposition; i.e.,
the relations $(\dagger)$, and therefore the earlier discussion
applies. First, we determine that the Cuntz-Krieger part of $\{
T_e \}_{e \in \E(G)}$ is trivial. Indeed,
\begin{eqnarray*}
\H_c &=& \bigwedge_{k \in \bbN} \, \bigvee_{|u | = k} T_{u}(\H_G)
\\
&=& \bigwedge_{k \in \bbN} \,\bigvee_{|u | = k} Y^{1/2}S_u (\H )        \\
     &=& Y^{1/2}\Big( \bigwedge_{k \in \bbN} \, \bigvee_{|u | = k} S_u (\H) \Big)
     = \{0 \}.
\end{eqnarray*}
Thus $\{ T_e \}_{e \in \E(G)}$ is pure. Let $\td{\W}$ be its
wandering subspace and let $\td{\W}_x \equiv P_x \td{\W} \subseteq
\td{\W}$. Then by the Wold decomposition there exists a unitary
operator
\[
V: \H \longrightarrow \sum_{x \in \V(G)}\oplus(Q_x\H_G)^{(\beta_x)}
\]
such that for any $e \in \E(G) $,
\[
VT_e V^* = \sum_{x \in \V(G)}\oplus L_{e}^{(\beta_x)}\Bigr| _{
(Q_x\H_G)^{(\beta_x)}}  ,
\]
where, $\beta_x = \dim \td{\W}_x$, $x \in \V(G)$. Similarly, there
exists a unitary
\[
U: \H \longrightarrow \sum_{x \in \V(G)}\oplus(Q_x\H_G)^{(\alpha_x)}
\]
such that for any $e \in \E(G) $,
\[
US_e U^* = \sum_{x \in \V(G)}\oplus L_{e}^{(\alpha_x)}\Bigr| _{
(Q_x\H_G)^{(\alpha_x)}}  .
\]
In order to establish the joint unitary equivalence between the
families $\{ T_e \}_{e \in \E(G)}$ and $\{ S_e \}_{e \in \E(G)}$,
we will prove that $VT_e V^* = U^* S_e U$ , $e \in \E(G)$. From
what we have proved so far, it follows that we need only to to
verify $\beta_x =\alpha_x$, for all $x \in \V(G)$. Fix an $x \in
\V(G)$ and let $\xi \in \td{\W}_x$. Then for all $ e \in \E(G)$
and $ \zeta \in \H_G$, we have
\[
( Y^{1/2}\xi \, , \, S_e \zeta )  =  ( \xi \, , \, Y^{1/2}S_e
\zeta ) =  ( \xi \, , \, T_e Y^{1/2} \zeta ) = 0.
\]
Thus by the invertibility of $Y$ we have  $Y^{1 / 2} \td{\W}  =
\W$. But $Y^{1/2} T_e^* T_e = S_e^* S_e Y^{1/2}$ for $e\in\E(G)$
since the same is true for $Y$. Hence, as $G$ has no sinks, every
$P_x$ is equal to some initial projection $S_e^*S_e =T_e^*T_e$ and
it follows that $Y^{1/2} \td{\W}_x = \W_x$ for $x \in\V(G)$ since
$Y^{1/2}$ is invertible. Therefore, the identity $\beta_x =
\alpha_x$ holds for all $x\in\V(G)$, as required.

We have shown that  $VT_e V^* = US_e U^*$ for $e\in\E(G)$. In
particular, this implies that  $U^*V Y^{1/2}S_e = S_e U^* V
Y^{1/2}$, for all $e \in \E(G)$, and the operator $U^* V Y^{1/2}$
commutes with $\{ S_e\}_{e \in \E(G)}$. Moreover, notice that for
$e\in\E(G)$,
\begin{eqnarray*}
S_e^* S_e U^*V Y^{1/2} &=&  S_e^*  U^*VT_e V^* V Y^{1/2} \\
&=& U^*V T_e^* V^* V T_e Y^{1/2} =  U^*V Y^{1/2} T_e^* T_e,
\end{eqnarray*}
and it follows that $U^* V Y^{1/2}$ commutes with
$\{P_x\}_{x\in\V(G)}$ as well. The desired operator for $(i)$ is
$A = U^* V Y^{1/2}$. \bx

Using the fact that $\Rg^{\prime} = \Lg$ we obtain the following.

\begin{cor}  \label{factorprop}
Let $G$ be a countable directed graph with no sinks. If $ Y \in
B(\H )$ is a positive invertible operator then the following
assertions are equivalent:
\begin{itemize}
\item[(i)] $Y = A^{*}A$ for some $A \in \Lg$. \item[(ii)]
$R_{e}^{*}YR_f = YR_{e}^{*} R_f = R_{e}^{*} R_f Y$ for all $e,f
\in \E(G)$.
\end{itemize}
\end{cor}

\section{An Application To The Geometry Of The Unit Ball}

In \cite{DKP}, Davidson, Pitts and the first author proved that
every operator in the open unit ball of a w*-closed operator
algebra generated by $n$ isometries with orthogonal ranges, is the
average of isometries. This result applies in particular to the
noncommutative Toeplitz algebra $\fL_n$, thus extending a
classical result of Marshall \cite{Mar} from function theory to
the noncommutative setting. Using the Szego type factorization
theorem we proved earlier, we extend the result of \cite{DKP} to a
broad class of free semigroupoid algebras.

The techniques of \cite{DKP}
yield a more general result than the one quoted in the introduction.
In order to state it, we need the following.

\begin{defn}
An operator
algebra $\A$ satisfies the \textit{factorization property} if for
any $A \in \A$ with $\| A\|< 1$, the equation $X^* X = I - A^* A$ has a
solution in $X \in \A$.
\end{defn}

It is clear from Corollary \ref{factorprop} that any free
semigroupoid algebra for a graph with no sinks satisfies the
factorization property.

\begin{thm}  \label{thm:DKP}
\cite{DKP} Let $\A$ be a norm closed operator algebra satisfying
the following two properties:
\begin{itemize}
\item[(i)] $\A$ satisfies the factorization property.
\item[(ii)] $\A$ contains two isometries with orthogonal ranges.
\end{itemize}
If $\|A\| < 1-\frac{1}{k}$, $A \in \A$, $k \in \bbN$, then
$A$ is the average of $6k$ isometries.
\end{thm}

\begin{rem}
The above result is not explicitly stated as a theorem in \cite{DKP}
but its validity is ascertained on \cite[page 118]{DKP}.
\end{rem}

On the other hand, in \cite{KP2}  Power and the second author
identified a property for a directed graph $G$ that is equivalent
to the existence of two isometries in $\Lg$ with orthogonal
ranges.

\begin{defn} \label{APP}
A countable directed graph satisfies the \textit{uniform aperiodic
path property} if the saturation of each vertex $x \in \V(G)$
either contains two distinct loops or an infinite proper path.
\end{defn}

\begin{thm} \label{aperiodic}
\cite{KP2} Let $G$ be a countable directed graph. Then $G$
satisfies the uniform aperiodic path property if and only if $\Lg$
contains a pair of isometries with orthogonal ranges.
\end{thm}

In light of Theorems \ref{thm:DKP} and \ref{factorprop}, the previous result provides the last
step for establishing the main result of this section.

\begin{thm}  \label{main}
Let $G$ be a countable directed graph which satisfies the uniform
aperiodic path property and let $A \in \Lg$. If $\|A\| <
1-\frac{1}{k}$, $k \in \bbN$, then $A$ is the average of $6k$
isometries.
\end{thm}

\begin{rem}
What about other free semigroupoid algebras? Clearly, $H^{\infty}$
fails the aperiodic path property and yet  the convex combinations
of the isometries cover the open unit ball. On the other hand, not
all free semigroupoid algebras satisfy that property. Indeed, if
$G$ is a graph with two vertices $x , y$ and an edge $e = xey$,
then $\Lg$ is the direct sum of the $2 \times 2$ upper triangular
matrices with $\bbC$. The convex hull of the isometries in this
algebra consists of the diagonal matrices.
\end{rem}

\section{Von Neumann Inequality}\label{S:von}

Given a row contraction $T = (T_1, \ldots,T_n)$ of (non-zero)
operators acting on a Hilbert space $\H$, consider a family of
mutually orthogonal projections $\P = \{ P_x \}_{x\in\J}$ on $\H$
which sum to the identity operator and  {\it stabilize} $T$ in the
following sense:
\begin{eqnarray}\label{stabilize}
T_i P_x ,\,\, P_x T_i \in \{ T_i, 0 \} \qfor 1\leq i \leq n \qand
x\in\J.
\end{eqnarray}
Observe that these relations determine a directed graph $G$ with
vertex set $\V(G) \equiv \J$ and $n$ directed edges $e_1, \ldots
e_n$, where $r(e_i)$ and $s(e_i)$ are the unique vertices with
$P_{r(e_i)} T_i P_{s(e_i)} = T_i$. When there is no confusion, we
shall write $T = (T_e )_{e\in\E(G)}$ if an ordering of the $T_i$
has been induced by a fixed projection set $\P$.

Let $T = \rowt$ and $\P = \{P_x\}_{x\in\V(G)}$  satisfy
(\ref{stabilize}) on $\H$ for a graph $G$ with no sinks. Then it
was shown in \cite{JK2}, and previously in a more abstract form
\cite{MS2} for the general setting of tensor algebras over
$\ca$-correspondences, that there are partial isometries $\{ S_e :
e\in\E(G)\}$ on a space $\K\supseteq\H$ such that
\begin{itemize}
\item[$(a)$] $S = \rows$ satisfy $(\dagger)$. \item[$(b)$] $\H$
reduces each $S_e^* S_e$, $e \in \E(G)$, and $\{S_e^*S_e|_\H\} =
\P$. \item[$(c)$] $\H$ is invariant for each $S_e^*$ with
$S_e^*|_{\H} = T_e^*$, $e\in\E(G)$. \item[$(d)$]
 $\K =  \bigvee_{w\in\fgeeplus} w(S) \H $.
\end{itemize}
Such a {\it minimal partially isometric dilation} of $T$ is unique
up to a unitary equivalence which fixes $\H$. The vertex
multiplicities $\alpha_x$ in such a dilation are determined by
(\ref{vertexmult}) and, in fact, it was shown in \cite{JK2} that
they may be computed from $T, \P$ by
\begin{eqnarray}\label{vertexmult2}
\alpha_x = \rank \Big( P_x \big( I_\H - \sum_e T_e T_e^* \big)
\Big) \qfor x\in\V(G).
\end{eqnarray}

For a fixed directed graph $G$ we shall define the supremum norm
$||p||_\infty$ of a polynomial $p\in\P_G^+$ by $||p||_\infty
\equiv ||p(L_G)||$, where $L_G = \rowl$ is the canonical pure
model (\ref{ledefn}) for $G$. This definition is in line with
other settings, where the norm $||p||_\infty$ is determined by the
pure model in this manner. As the simplest example, recall that
the norm of a polynomial $p$ on the unit complex disk satisfies
$||p||_\infty =||p(U_+)||$, where $U_+$ is the unilateral forward
shift \cite{SzF}.

\begin{thm}\label{vonthm}
Let $G$ be a countable directed graph with no sinks. Let $T=
(T_e)_{e\in \E(G)}$ be a row contraction on a Hilbert space $\H$
which satisfies (\ref{stabilize}) for $\{P_x:x\in\V(G)\}$. Then
\begin{eqnarray}\label{ineq}
|| p(T)  || \leq || p ||_\infty \qforal p\in \P_G^+.
\end{eqnarray}
\end{thm}

\Prf Let $S=(S_e)_{e\in\E(G)}$ be the minimal partially  isometric
dilation of $T$ associated with the projections $\{P_x\}$. Let $p$
belong to $ \P_G^+$. Then $p(T) = P_\H \, p(V)|_\H$ by virtue of
$(a)$, $(b)$ and $(c)$ above, and hence the Wold decomposition
implies that
\[
|| p(T)  || \leq || p(S)  || = \max \{ || p(V)  ||, || p(W)  ||\},
\]
where $V=(V_e)_{e\in\E(G)}$ satisfies $(\ddagger)$ and
$W=(W_e)_{e\in\E(G)}$ is pure and satisfies $(\dagger)$.

The decomposition (\ref{wold}) of the pure part $W$ of the minimal
dilation $S$ clearly implies that
\[
|| p(W)  || \leq || p(L_G)  || \equiv || p ||_\infty.
\]
Thus, we shall finish the proof by showing that $|| p(S)  || \leq
|| p(L_G) ||$ for {\it every} row contraction
$S=(S_e)_{e\in\E(G)}$ which satisfies $(\ddagger)$ for $G$.

To see this, let $S=(S_e)_{e\in\E(G)}$ satisfy $(\ddagger)$ and
let $0 < r < 1$. Then $rS = (rS_e)_{e\in\E(G)}$ is a pure row
contraction, and hence the minimal partially isometric dilation of
$rS$ with respect to $\{S_e^*S_e : e\in\E(G)\}$ is given by the
form (\ref{wold}) with $V_e \equiv 0$. We claim that the vertex
multiplicities for this pure dilation
$S'=(S_e^\prime)_{e\in\E(G)}$ satisfy $\alpha_x \geq 1$ for
$x\in\V(G)$. Indeed, from (\ref{vertexmult2})  these
multiplicities are determined by
\[
\alpha_x = \rank \Big(P_x - \sum_{e=xe} (rS_e) (rS_e)^* \Big)
\qfor x\in\V(G).
\]
Thus if $\alpha_x =0$, then the $(\ddagger)$ relations give
\[
\sum_{e=xe} S_e S_e^* = P_x = r^2 \sum_{e=xe} S_e S_e^*,
\]
which is clearly a contradiction since the $P_x$ are assumed to be
non-zero in the $(\dagger)$ relations.

As $\alpha_x \geq 1$ for $x\in\V(G)$, we have
\[
|| p(rS) || \leq || p(S') || =|| p(L_G)  ||.
\]
Since $0 < r < 1$ was arbitrary, it follows that $|| p(S)  || \leq
|| p(L_G)  ||$, as required. \bx

\begin{rem}
In the case of single vertex graphs with $n$ edges, the above
theorem recovers the classical von Neumann inequality
\cite{vN,SzF} (when $n=1$) and Popescu's version \cite{Popvn}
(when $2\leq n \leq \infty$). We also mention that the conclusion
of Theorem~5.1 can be proved through an application of
Theorem~3.10 from \cite{MS2}. However, we note that our direct
proof via the Wold decomposition is novel and provides a new
perspective on the problem.
\end{rem}

\subsection{Partially Ordered Sets of Directed
Graphs}\label{sS:partial} 

An alternative perspective on  Theorem~\ref{vonthm} is the
following: Given a countable directed graph $G$ with no sinks, the
theorem shows there is a corresponding von Neumann inequality with
estimate given by (for this discussion) $||p||_{G,\infty} \equiv
||p(L_G)||$ for $p\in\P_G^+$, recalling that  $L_G =
(L_e)_{e\in\E(G)}$. On the other hand, every directed graph
determines a partially ordered set of directed graphs through its
deformations. A directed graph $G_2$ is a {\it deformation} of
$G_1$ if $G_2$ is obtained from $G_1$ by identifying certain
vertices in $G_1$. The partial order is defined by: $G_1\leq G_2$
if and only if $G_2$ is a deformation of $G_1$. (This perspective
was also discussed in \cite{JK2} in the context of dilation
theory, where it was shown that every row contraction $T$
generates such a partially ordered set of directed graphs through
its family of minimal partially isometric dilations.) We observe
below that the von Neumann inequality estimates respect the
ordering in such a partially ordered set of directed graphs.

Let $C_n$ be the directed graph with a single vertex and $n$
directed loop edges. The Sz.-Nagy and Popescu estimates are given
by $||p(L_{C_n})||$, where $p$  is allowed to be any polynomial in
$n$ noncommuting variables when $n\geq 2$. Let $G$ be a countable
directed graph with $n$ edges. Then $L_G = \rowl$ is a pure row
contraction (in the $n$-tuple sense \cite{Pop_diln}) and hence its
(unique) minimal isometric dilation \cite{Bun,Fra1,Pop_diln} is
given by a multiple of $L_{C_n} = (L_1,\ldots ,L_n)$, where the
$L_i$ are the creation isometries on the full $n$-variable Fock
space. Thus, given $L_G = \rowl$ and such a polynomial $p$, we
shall write $p(L_G)$ for the evaluation of $p$ at the $L_e$ which
arises from the ordering of the $L_e$ induced by the minimal
isometric dilation of $L_G$.

\begin{cor}\label{poset}
If $G_1 \leq G_2$ are directed graphs with no sinks and $n$
directed edges, then for every polynomial $p$ in $n$ noncommuting
variables
\[
||p||_{G_1,\infty} \leq ||p||_{G_2,\infty} \leq
||p||_{C_n,\infty}.
\]
\end{cor}

\Prf The relation $G_1 \leq G_2$ can be seen to induce a
$G_2$-ordering on $L_{G_1}$, where the corresponding projections
which stabilize $L_{G_1}$ are given by the sums of the initial
projections for $L_e$, $e\in\E(G_1)$, which are naturally induced
by the deformation. Hence the first inequality is a consequence of
Theorem~\ref{vonthm}.  The last inequality follows from Popescu's
von Neumann inequality since $L_{G_1}$ and $L_{G_2}$ are both row
contractions and $n$-tuples. \bx

It is not hard to see that $||p||_{G,\infty}$ gives a sharper
estimate than $||p(L_{C_n})||$ in general when $n\geq 2$. Consider
the following simple illustration of this fact.

\begin{eg}
Let $T_1,T_2$ be contractions on a Hilbert space $\H$ and define
operators on $\H^{(2)}$ by
\[
V_1= \left( \begin{matrix}
0 & 0 \\
T_1 & 0
\end{matrix}\right)  \quad {\rm and} \quad
V_2= \left(\begin{matrix}
0 & T_2 \\
0 & 0
\end{matrix}\right).
\]
Then $V=(V_1,V_2)$ is a row contraction which is stabilized by $\P
= \{ P_1, P_2\}$ where $P_i \equiv P_\H$ are the projections of
$H^{(2)}$ onto its coordinate spaces. Let $G$ be the graph with
two vertices $x,y$ and edges $e=yex$, $f=xfy$. Then $L_G =
(L_e,L_f)$ in this case. Observe that $C_2$ is a deformation (the
only one) of $G$.  Let $L_{C_2} =(L_1,L_2)$, where $L_1,L_2$ are
the creation isometries on the full 2-variable Fock space. Then
$L_e + L_f$ is an isometry and $L_1 + L_2 = \sqrt{2}L$ where $L$
is an isometry, and hence
\[
|| V_1 + V_2 || \leq || L_e + L_f || =1 < \sqrt{2} = || L_1 + L_2
||.
\]
Further, the $L_e+L_f$ estimate is best possible in this example
since $||V_1 + V_2|| =1$ if $||T_i||=1$ for $i=1$ or $2$.
\end{eg}

\section{Functional Calculus}\label{S:func}

The von Neumann inequality leads to a natural functional calculus
determined by quiver algebras and, in the case of pure row
contractions, by free semigroupoid algebras.  Observe that when
$T= (T_e)_{e\in \E(G)}$ is a row contraction with $G$-ordering
induced by projections $\P=\{P_x\}$ which satisfy
(\ref{stabilize}) for $T$, we may consider all polynomials $\{
p(T) : p \in \P_G^+ \}$ in the generators $T_e,P_x$. The following
result is an immediate consequence of Theorem~\ref{vonthm} and may
also be obtained as a special case of Theorem~3.10 from
\cite{MS2}.

\begin{thm}\label{fcalc}
Let $G$ be a countable directed graph with no sinks. Let $T=
(T_e)_{e\in \E(G)}$ be a row contraction on $\H$ which satisfies
(\ref{stabilize}) for  $\{P_x^\prime:x\in\V(G)\}$. Then there is a
contractive homomorphism
\[
\Psi : \A_G \longrightarrow \bofh
\]
defined by $\Psi(L_e) = T_e$ and $\Psi(P_x) = P_x^\prime$ for
$e\in\E(G)$ and $x\in\V(G)$.
\end{thm}

When the row contraction satisfies an extra condition, the
functional calculus can be extended to $\flgee$.

\begin{thm}\label{fcalcwot}
Let $G$ be a countable directed graph with no sinks. Let $T=
(T_e)_{e\in \E(G)}$ be a row contraction on $\H$ which satisfies
(\ref{stabilize}) for  $\{P_x^\prime:x\in\V(G)\}$, and suppose
that
\begin{eqnarray}\label{purecond}
\lim_{d\rightarrow \infty} \Big( \sum_{w\in\fgeeplus_d} || w(T)^*
\xi || \Big) = 0 \qfor \xi\in\H.
\end{eqnarray}
Then there is a contractive homomorphism
\[
\Psi : \fL_G \longrightarrow \bofh
\]
defined by $\Psi(L_e) = T_e$ and $\Psi(P_x) = P_x^\prime$ for
$e\in\E(G)$ and $x\in\V(G)$.
\end{thm}

\Prf Condition (\ref{purecond}) means precisely that the minimal
partially isometric dilation $V = (V_e)_{e\in\E(G)}$ of $T$ with
respect to $\{ P_x^\prime \}$ is pure \cite{JK2}. By the Wold
decomposition for $V$, we may define operators $f(V)$ for all
$f\in\flgee$, where  $f(V) = \wot\!-\!\lim_k p_k(V)$ for some
sequence $p_k \in \P_G^+$. From the properties of the minimal
dilation we have $p_k(T) = P_\H \, p_k(V)|_\H$. Thus, $p_k(T)$
converges to an operator $f(T) \equiv \wot\!-\!\lim_k
p_k(T)\in\bofh$. Finally, by the von Neumann inequality and a
standard approximation argument we have $||f(T)|| \leq ||f(V)||
\leq ||f(L)||$, and the result follows. \bx

\begin{rem}
It should be possible to prove an $\flgee$ functional calculus for
more general row contractions. For instance, as in
\cite{SzF,BV1,Pop2} this could be accomplished by first
identifying `characteristic functions' for an appropriate notion
of {\it completely non-unitary} row contractions here.
\end{rem}




\section{A Generalization of Coburn's Theorem}\label{S:Coburn}

In \cite{Cob} Coburn showed that the $\ca$-algebra generated by a
non-unitary isometry $V$ is isomorphic to the Toeplitz
$\ca$-algebra, via a map that sends $V$ to the forward shift. This
result was generalized by Popescu to Cuntz-Toeplitz $\ca$-algebras
\cite{Pop2}. In this section we consider the generalization of
Coburn's Theorem to $\ca$-algebras generated by partial isometries
associated with more general directed graphs.

The proof of Coburn's theorem depends on the universality of the
algebra of continuous functions on the circle as the $\ca$-algebra
of a graph $G$ consisting of one vertex and one loop-edge. The
following result of Szymanski \cite{Szy} generalizes this fact to
arbitrary graph $\ca$-algebras.

\begin{thm} \label{Szy}
\cite{Szy} Let $G$ be a countable directed graph with no sources
and $\rho : \ca(G) \longrightarrow \A$ be a $*$-homomorphism into
a $\ca$-algebra $\A$. Then $\rho$ is injective if and only if the
following two conditions are satisfied:
\begin{itemize}
\item[(i)] $\rho (P_x) \neq 0$ for each vertex $x \in \V(G)$.
\item[(ii)] For each vertex simple loop $w \in \fgeeplus$ without
entrances, we have $\bbT \subseteq \sigma( \rho(S_w))$.
\end{itemize}
\end{thm}

\begin{cor} \label{Kquotient}
Let $G$ be a countable directed graph with no  sources and no
sinks. Let $\ca (\A_G)$ be the $\ca$-algebra generated by the
quiver algebra $\A_G$, let $\K$ denote the set of compact
operators on $\H_G$ and let $ \ca (\A_G) \slash \, \K$
denote the image of $ \ca (\A_G)$ under the Calkin map.
Then there exists a $*$-isomorphism
\[
\rho : \ca (G) \rightarrow \ca (\A_G) \slash \, \K
\]
such that $\rho( S_e ) = L_e + \K$, for all $e \in \E(G)$.
\end{cor}

\Prf Since $G$ has no sources, the projections $L_x$, $ x \in
\V(G)$, are infinite dimensional and thus $L_x + \K \neq 0$.
Furthermore,
\[
\sum_{r(e) = x} L_e L_{e}^{*} = L_x - \xi_x \otimes \xi_x
\]
and so the collections $\{ L_x + \K \}_{x \in \V(G)}$ and $\{ L_e
+ \K \}_{e \in \E(G)}$ satisfy the relations $(\ddagger)$. Hence
there exists a $*$-homomorphism $\rho : \ca (G) \longrightarrow
\ca(\A_G) \slash \, \K$ with $\rho( P_x ) = L_x + \K$ and $\rho(
S_e ) = L_e + \K$ for all $x \in \V(G)$ and $e \in \E(G)$. We now
use Szymanski's result to show that $\rho$ is injective.

Towards this end, notice that we have already verified condition
(i) in Theorem~\ref{Szy}. To verify the other condition let $w \in
\fgeeplus$ be any loop. Consider the direct sum decomposition
$\H_G = L_x (\H_G ) \oplus L_x (\H_G )^{\perp}$ and notice that
with respect to this decomposition, $L_w$ admits the $2 \times 2$
matrix form
\[
\left(
\begin{matrix}

L & 0 \\
0 & 0

\end{matrix}
\right) ,
\]
with $L$ unitarily equivalent to some ampliation
$S^{(\alpha )}$ of the forward shift $S$. It is well known that if
$\lambda \in \bbT$, then the range of $S - \lambda I$ is not
closed. Therefore the same is true for  $S^{(\alpha )}- \lambda I$
and so $\lambda$ belongs to the essential spectrum of $S^{(\alpha
)}$; i.e., $\lambda \in \sigma( \rho(S_w) )$, as desired. \bx

With a standard `diagram chasing' argument,
Corollary~\ref{Kquotient} can be used to show the essential norm
coincides with the operator norm for any element of $\A_G$.
However, we can prove this directly for a larger algebra, the free
semigroupoid algebra $\flgee$, hence we include a proof of this
fact.

\begin{prop} \label{essnorm}
The following assertions are equivalent for a countable directed
graph $G$:
\begin{itemize}
\item[$(i)$] $G$ has no sources.
\item[$(ii)$] $\| A \|
=\|A\|_{e}$ for all $A\in \flgee$.
\end{itemize}
\end{prop}

\Prf If $x$ is a source vertex in $G$, then $P_x$ is a rank one
projection inside $\A_G$ (and $\flgee$). Conversely, suppose $G$
has no sources. Let $d\geq 1$ be a positive integer. For all
$x\in\V(G)$ we may choose a path $w_x\in\fgeeplus$ of length $d$
such that $r(w_x) =x$. As the $w_x$ are distinct paths of the same
length, notice that the partial isometries $ R_{w_x}$,
$x\in\V(G)$, have pairwise orthogonal ranges. Thus, we may define
an isometry $R^{(d)} = \sum_{x\in\V(G)} R_{w_x}$ in the commutant
$\flgee^\prime = \frgee$ (where the sum converges $\wot$ in the
infinite vertex case);
\[
(R^{(d)})^* R^{(d)} = \sum_{x,y} R_{w_x}^* R_{w_y} = \sum_x
R_{w_x}^* R_{w_x} = \sum_x Q_x = I.
\]

Hence, if $A$ belongs to $\flgee$ we have
\[
||AR^{(d)} \xi|| = || R^{(d)} A \xi|| = ||A\xi|| \qfor \xi\in\H_G
\qand d\geq 1.
\]
On the other hand, $\{ R^{(d)}\xi\}_{d\geq 1}$ converges weakly to
zero for any choice of vector $\xi\in\H_G$, and thus $||A|| =
||A||_e$ as claimed. \bx

\begin{cor}
If $G$ be a directed graph with no sources, then the compact perturbations
$\flgee + \K$ are norm closed.
\end{cor}

The next result shows that $\ca(\A_G)$ is the universal
$\ca$-algebra for the relations $(\dagger)$ and it is
known to the specialists (compare with \cite{P, H}). Nevertheless, we feel that the
short proof
presented here is new and we therefore include it.

\begin{thm}\label{daggeruniversal}
Let $G$ be a countable directed graph with no sources and no
sinks. Let $\{ P_{x}^{\prime} \}_{x \in \V(G)}$ and $\{
S_{e}^{\prime} \}_{e \in \E(G)}$ be families of projections and
partial isometries acting on a Hilbert space $\H$ which satisfy
$(\dagger)$. Then there exists a $*$-epimorphism
\[
\tau : \ca (\A_G ) \longrightarrow \ca (  \{ S_{e}^{\prime} \}_{e
\in \E(G)} )
\]
such that $\tau(L_e) = S_{e}^{\prime}$, for all $e \in \E(G)$.
\end{thm}

\Prf From the discussion in Section~\ref{S:wold}, it is clear that
the Wold decomposition  (which was proved  for $G$ with no sinks)
implies the existence of $\ca$-homomorphisms
\[
\phi_1 : \ca(G) \longrightarrow B (\H_1)
\]
and
\[
\phi_2 : \ca(\A_G) \longrightarrow B (\H_2)
\]
such that $\H = \H_1 \oplus \H_2$ and
\[
S^{\prime}_{u} \simeq (\phi_1 \oplus \phi_2)(S_u \oplus L_u )
\qfor u \in \fgeeplus.
\]
Consider the diagram
\[
\begin{CD}
\ca(\A_G)@> \pi>>  \ca(\A_G)\slash \K @> \rho>> \ca(G)  ,
\end{CD}
\]
where $\rho $ is the $*$-isomorphism identified by
Corollary~\ref{Kquotient}. If $\psi \equiv \rho \circ \pi$, then
the desired map here is
\[
\tau = (\phi_1 \oplus \phi_2) \circ (\psi \oplus id) : \ca(\A_G)
\longrightarrow \ca ( \{ S_{e}^{\prime} \}_{e \in \E(G)} ).
\]
Tracing the definitions of these maps shows that $\tau (L_u) =
S^{\prime}_{u}$ for all $ u \in \fgeeplus$. \bx

The previous theorem leads to a natural $\ca$-algebra extension of
the von Neumann inequality Theorem~\ref{vonthm} (see also
\cite{Pop2} for the case of isometries with orthogonal ranges).
Though note that the proof of this result relies on Szymanski's
Theorem, whereas our proof of Theorem~\ref{vonthm} is direct. Let
$\P_G^0$ be the set of polynomials generated by all indeterminant
monomials of the form $v\ol{w}$ with $v,w\in\fgeeplus$ and where
$x=\ol{x}$ for $x\in\V(G)$. Given $T = \rowt$ and $p\in\P_G^0$, we
write $p(T,T^*)$ for the polynomial $p$ evaluated at the
generators $T_e$ (for $e\in\P_G^0$) and $T_e^*$ (for
$\ol{e}\in\P_G^0$) and $P_x$ (for $x=\ol{x} \in\P_G^0\cap\V(G)$).

\begin{cor}
Let $G$ be a countable directed graph with no sources and no
sinks. Let $T = \rowt$ be a row contraction on $\H$ stabilized as
in (\ref{stabilize}) by a given family of projections. Then
\begin{eqnarray}\label{vnextended}
|| p(T,T^*) || \leq || p(L_G , L_G^*) || \qfor p\in\P_G^0.
\end{eqnarray}
\end{cor}

\Prf Let $S = \rows$ be the minimal partially isometric dilation
of $T$ with respect to the given projection family. From the
basic properties of the minimal dilation, it follows that
\[
p(T,T^*) = P_\H \, p(S,S^*)|_\H \qfor p\in\P_G^0.
\]
Thus, Theorem~\ref{daggeruniversal} gives $|| p(T,T^*) || \leq
||p(S,S^*)|| \leq || p(L_G , L_G^*) ||$ for all $p\in\P_G^0$. \bx

The next theorem is the promised generalization of Coburn's
Theorem. Given a collection $\{S_{e}^{\prime}\}_{e \in \E(G)}$ of
isometries, we isolate a condition which is the appropriate
extension of an isometry being non-unitary here (condition
(\ref{properisom}) below). We show that any collection of partial
isometries satisfying this condition generates a $\ca$-algebra
isomorphic to $\ca(\A_G)$.

\begin{thm}   \label{Coburn}
Let $G$ be a countable directed graph with no  sources and no
sinks and let $\ca (\A_G)$ be the $\ca$-algebra generated by the
quiver algebra $\A_G$. Let $\{P_{x}^{\prime}\}_{x \in \V(G)}$ and
$\{S_{e}^{\prime}\}_{e \in \E(G)}$ be families of projections and
partial isometries respectively, acting on a Hilbert space $\H$
and satisfying $(\dagger)$. If,
\begin{equation} \label{properisom}
\sum_{r(e) = x} S_{e}^{\prime} (S_{e}^{\prime})^{ *} \neq
P_{x}^{\prime}     \qfor x \in \V(G),
\end{equation}
then there exists an injective $*$-homorphism $\phi: \ca (\A_G )
\rightarrow B(\H)$ such that $\phi( L_e ) = S_{e}^{\prime}$ for
all $e \in \E(G)$.

In particular, the $\ca$-algebra $\ca(\{S_{e}^{\prime}\}_{e \in
\E(G)})$, generated by the collection $\{S_{e}^{\prime}\}_{e \in
\E(G)}$, is isomorphic to $\ca (\A_G)$.
\end{thm}

\Prf We utilize the Wold decomposition for the collection
$\{S_{e}^{\prime}\}_{e \in \E(G)}$. Recall the existence of
reducing subspaces $\H_c$ and $\H_p$ for $\{S_{e}^{\prime}\}_{e
\in \E(G)}$ such that $\H = \H_p \oplus \H_c$ and the restriction
of $\{S_{e}^{\prime}\}_{e \in \E(G)}$ on $\H_p$ (respectively
$\H_c$) is pure (respectively fully coisometric). Moreover, if
\[
\W \equiv \ran \Big( I\, -\, \sum_{e} S_{e}^{\prime}
(S_{e}^{\prime})^{ *}\Big)
\]
then $\W$ is wandering for $\{S_{e}^{\prime}\}_{e \in \E(G)}$ and
the smallest invariant subspace for $\{S_{e}^{\prime}\}_{e \in
\E(G)}$ containing $\W$ equals $\H_p$. In addition, the
restrictions $S_{e}^{\prime}|_{\H_p}$, $e \in \E(G)$, are jointly
unitarily equivalent to
\[
S_{e}^{\prime}|_{\H_p}\simeq \sum_{x \in \V(G)}\, \oplus
L_{e}^{(\alpha_x )}|_{\, ( Q_x \H_G )^{(\alpha_x )}},
\]
where,
\[
\alpha_x = \rank\left( P_{x}^{'} - \sum_{r(e) = x} S_{e}^{\prime}
(S_{e}^{\prime})^{ *}\right).
\]

Let $\W_x = P_{x}^{\prime} \W \neq 0$ and let $\zeta_x \in \W_x$
be a unit vector. Notice that the subspaces
\[
\M_x   \equiv  \spn \{ S_{w}\zeta_x : w \in \fgeeplus \}
\]
are reducing for $\{S_{e}^{\prime}\}_{e \in \E(G)}$ and the
restrictions $S_{e}^{\prime} |_{M_x}$, $e \in \E(G)$, are jointly
unitarily equivalent to $L_{e}|_{\, Q_{x}}$, $e \in \E(G)$,
respectively. Therefore, the subspace $\M =\sum \oplus_{x \in
\V(G)}\M_x$ is reducing for $\{S_{e}^{\prime}\}_{e \in \E(G)}$ and
the restrictions $S_{e}^{\prime}|_{\, \M_x}$, $e \in \E(G)$, are
jointly unitarily equivalent to $L_e$, $e \in \E(G)$. We have
therefore established the existence of an epimorphism
\[
\psi: \ca(\{S_{e}^{\prime}\}_{e \in \E(G)})  \longrightarrow
\ca(\A_G)
\]
mapping generators to generators. It remains to show that $\psi$
is injective.

Consider the diagram
\[
\begin{CD}
\ca(\A_G)@> \pi>>  \ca(\A_G)\slash \K @> \rho>> \ca(G)  @>  \tau
>> \ca(\{S_{e}^{\prime}|_{\, \H_{c}}\}_{e \in \E(G)}),
\end{CD}
\]
where $\pi$ is the natural quotient map, $\rho$ is as in
Lemma~\ref{Kquotient} and $\tau$ exists because $\ca(G)$ is
universal. Let
\[
\phi_1 \equiv \tau \circ \rho \circ \pi : \ca(\A_G)
\longrightarrow \ca(\{S_{e}^{\prime}|_{\, \H_{c}}\}_{e \in \E(G)})
\]
and note that $\phi_1 (L_e) = S_{e}^{'}|_{\H_c}$, for all $e \in
\E(G)$. Consider now the map $\phi :\ca(\A_G) \longrightarrow
\ca(\{S_{e}^{\prime}\}_{e \in \E(G)})$, defined as
\[
p( L_e , L_{e}^{*} ) \longmapsto p( S_{e}^{\prime} ,
(S_{e}^{\prime})^{*} )|_{\H_p} \oplus \phi_1( p( L_e , L_{e}^{*}
)),
\]
for any non-commutative polynomial $p$. The Wold decomposition
shows that $\phi$ is a $*$-homomorphism and the earlier discussion
shows that $\phi(L_e) = S_{e}^{\prime}$, for all $e \in \E(G)$.
Therefore, $\phi \circ \psi = id$ and so $\psi$ is injective,
which proves the theorem. \bx

\begin{rem}
It was pointed to us by the referee that Theorem \ref{Coburn} has
been previously obtained by Fowler and Raeburn \cite[Theorem
4.1]{FR}. Their approach is quite different from ours since they
make an extensive use of Hilbert $\ca$-modules as well as their
own methods \cite{FR2} for analyzing semigroup crossed products.
Our proof is shorter and elementary; we only require Szymanski's
result, whose proof in  \cite{Szy} is essentially self-contained.
\end{rem}



\vspace{0.1in}

{\noindent}{\it Acknowledgements.} The first author was partially
supported by a research grant from ECU. The second author was
partially supported by an NSERC research grant and start up funds
from the University of Guelph. We would like to thank the referee
for a number of helpful comments.

\end{document}